\documentclass[12pt]{article}
\newtheorem{thm}{Theorem}[section]
\newtheorem{prop}[thm]{Proposition}
\newtheorem{cor}[thm]{Corollary}
\newtheorem{lemma}[thm]{Lemma}

\def\proof{{\sc Proof. }}

\usepackage{amsfonts}
\usepackage{amssymb}

\def\CP{\mathbb{C}{\rm P}} 
\def\Im{\mathop{\rm Im}}
\def\Re{\mathop{\rm Re}}
\def\R{\mathbb{R}} 
\def\C{\mathbb{C}}
\def\H{\mathbb{H}}
\def\Lc{{\cal L}}

\def\Sc{{\cal S}}   
\def\Box{\square}

\title{Rectification of circles and quaternions}
\author{V. Timorin
\thanks{Partially supported by RFBR 99-01-00245 and CRDF RM1-2086}}
\date{}

\begin{document}
\maketitle

{\small Consider a bundle of circles passing through 0 in 4-dimensional space. 
It is said to be {\em rectifiable} if there is a germ of diffeomorphism at 0 
that takes all circles from our bundle to straight lines. 
We will give a classification of all rectifiable bundles 
of circles containing sufficiently many circles in general 
position. This result is surprisingly different from those in dimensions 2 and 3 
(Khovanskii and Izadi) due to a connection with the quaternionic algebra.}

\section*{Introduction}

Throughout this paper, the word ``circle'' means a circle or a 
straight line. We are always assuming that the space $\R^n$ is equipped with
a fixed ``standard'' Euclidean inner product.  

A collection of curves in $\R^n$ passing through 0 is said to be a {\em simple
bundle of curves} if no two of them are tangent at 0.  
A simple bundle of curves is called {\em rectifiable} if there exists a germ of 
diffeomorphism in a neighborhood of the origin that sends all curves from this bundle to
straight lines. Rectifiable bundles of curves appear, for example, in Riemannian
geometry --- the set of geodesics passing through a given point is rectifiable.

A. G. Khovanskii  proved in \cite{Kh} that a rectifiable simple bundle
of more than 6 circles on plane necessarily pass through some 
point different from the origin. F. A. Izadi \cite{Iz} generalized 
Khovanskii's arguments to dimension 3. A rectifiable simple bundle of 
circles in $\R^3$ containing sufficiently many circles in general position
must pass through some other common point.    

In dimension 4, this is not true. The simplest counterexample is a family 
of circles that are obtained from straight lines by some complex 
projective transformation (with respect to some identification $\R^4=\C^2$
such that the multiplication by $i$ is an orthogonal operator).  

It turns out that in dimension 4 there is a large family of 
transformations that round lines (i.e., take them to circles). To
construct such a
family, fix a quaternionic structure on $\R^4$ compatible with the 
Euclidean structure. If $A$ and $B$ are some
affine maps, then the map $x\mapsto A(x)^{-1}B(x)$ rounds lines
(the multiplication and the inverse are in the sense of quaternions).
Such transformations will be called (left) {\em quaternionic fractional
transformations}.
Right quaternionic fractional transformations $AB^{-1}$ also round lines.
Any real projective, complex projective or quaternionic projective 
transformation is quaternionic fractional.  

In this paper, we will prove that a rectifiable simple bundle of circles
containing sufficiently many circles in general position is the image 
of a bundle of straight lines under some left or right quaternionic fractional 
transformation. 

In arbitrary dimension, we have a purely algebraic description of 
rectifiable simple bundles of circles. So the analytic problem 
of classification of such bundles is reduced to an algebraic problem.

The paper is organized as follows. In Section 1, for a simple rectifiable 
bundle of circles we establish an algebraic 
condition on the second derivative of a rectifying map. This condition
is formulated on the asymptotic cone $\{(x,x)=0\}\subseteq\C^2$ where $(\cdot,\cdot)$ is 
the complexification of the usual inner product. 
This provides a simple proof of Izadi's theorem \cite{Iz}. In Section 2, we show
that this algebraic condition is not only necessary but also sufficient in 
a sense. Thus we obtain an algebraic description of rectifiable simple bundles
of circles. In Section 3, we review some important properties of complex and 
quaternionic structures and relate them to the geometry of the asymptotic cone.
In Section 4, we define quaternionic fractional transformations and list some 
of their properties. Section 5 contains the main rectification result and 
some its geometrical consequences.

I am grateful to A. G. Khovanskii for useful discussions. 
      
\section{Rectifiable collections of circles}

The following result is true in dimensions 2 \cite{Kh} and 3 \cite{Iz}.

\begin{thm}
\label{rect23}
Consider a simple bundle of circles in $\R^2$ or $\R^3$ containing
sufficiently many circles in general position. If this bundle is rectifiable,
then all its circles pass through a common point different from the origin.  
\end{thm}

On plane, it is enough to take 7 circles. Theorem \ref{rect23} means, in 
particular, that if a generic family of circles can be rectified anyhow, then
it can be rectified by means of some inversion. As we will see later, this 
violates in dimension 4.

We need the following very simple lemma:

\begin{lemma}
\label{pr}
Consider a polynomial map $F:\R^n\to\R^n$ such that $F(x)$ is everywhere
proportional to $x$. Then $F(x)=G(x)x$ for some polynomial function
$G:\R^n\to\R$. If $F$ is homogeneous, then so is $G$.  
\end{lemma}

\proof Introduce a coordinate system $(x_0,\dots,x_{n-1})$.
Denote by $F_i$ the $i$-th component of $F$. Then the proportionality
condition reads as $x_iF_0-x_0F_i=0$. In particular, $F_0$ is
divisible by $x_0$. Denote the quotient by $G$. Then from our equation
we see that $F_i=Gx_i$. The last statement of the lemma is obvious. $\Box$

Extend the standard 
inner product $(\cdot,\cdot)$ from $\R^n$ to $\C^n$ by complex bilinearity.
The locus $(x,x)=0$ is called the {\em asymptotic cone}. Denote this 
cone by $C$. The asymptotic cone describes the behavior of circles
at infinity. Namely, any nondegenerate circle (not a line) is asymptotic to $C$. 

Let $\Phi:(\R^n,0)\to (\R^n,0)$ be a germ of diffeomorphism at 0 that sends several 
lines passing through the origin to circles. Suppose that the number of lines is big 
enough and that they are in general position. Denote this set of lines by $\Lc$. We 
can assume without loss of 
generality that $d_0\Phi=id$. To arrange this it is enough to compose $\Phi$ with some 
linear transformation (namely, the inverse of $d_0\Phi$) which certainly takes lines 
to lines. Let $\Phi=x+\Phi_2(x)+\cdots$ near 0 where $\Phi_2$ denotes the second order
terms. 

\begin{prop}
\label{cond}
The quadratic map $\Phi_2$ satisfies the following relations on the
asymptotic cone:
$$(\Phi_2(x),\Phi_2(x))=0,\quad (\Phi_2(x),x)=0.$$
\end{prop}

This proposition means that $\Phi_2$ preserves the asymptotic cone and
takes each vector $x\in C$ to a vector $y\in C$ such that $x$ and $y$ span
a subspace lying entirely in $C$. To give an informal explanation of this
result let us assume the following: 
\begin{itemize}
\item
The diffeomorphism $\Phi$ takes germs of {\em all} lines passing through
0 to germs of circles.
\item
Our diffeomorphism can be extended to a neighborhood of the
origin in $\C^n$ as a local holomorphic map.
\end{itemize}

Then $\Phi$ takes germs of complex lines to germs of {\em complex circles}.
By a complex circle we mean an algebraic curve that is given by equations
of the same form as for real circles but in complex variables and
with complex coefficients. In particular, we can define a complex
circle as the intersection of a complex 2-dimensional plane
in $\C^n$ and a complex sphere $(x-a,x-a)=R^2$ where $a\in\C^n$ and
$R\in\C$. Thus any complex circle is a plane curve of degree 2.

Take a complex line $L$ from the cone $C$. Then $\Phi(L)$ is
a complex circle. We know that this circle is tangent to $L$ at 0
and asymptotic to $C$ at infinity. Denote by $M$ the plane
where $\Phi(L)$ lies. Then either $M$ is contained
in $C$ or $M\cap C$ is a pair of intersecting lines in $M$.
In the latter case $\Phi(L)$ must coincide with one of these lines.
Indeed, $\Phi(L)$ intersects both lines at the origin and is asymptotic
to one of them. But a plane curve of degree 2 cannot intersect its own
asymptotic line. Note that $L$ is clearly in $M\cap C$, so
$\Phi(L)=L$. 

In any case, $L$ and $\Phi(L)$ span a vector subspace lying entirely in $C$.
Hence $\Phi_2(L)$ lies in this subspace. From this the proposition follows.

The above argument can be extended to a rigorous proof but, to
give a shorter proof, we will use another idea. 

\proof Make the inversion $I$ with respect to the origin and consider the
composition $I\circ\Phi$. 
The diffeomorphism $\Phi$ takes a line from $\Lc$ to a tangent circle (due
to the condition $d_0\Phi=id$) and $I$ sends circles or lines tangent at 0
to parallel lines. Therefore, $I\circ\Phi$ maps each line from $\Lc$ to a
parallel line.

Consider the Taylor series for $\Phi$ at the origin:
$$\Phi(x)=x+\Phi_2(x)+\Phi_3(x)+\cdots,$$
where $\Phi_k(x)$ denotes the order-$k$ terms. Fix some nonzero vector $x$
that spans a line from $\Lc$.
This line can be parameterized as $\{xt\}$ where $t$ is a parameter.
Hence $I\circ\Phi(xt)$ runs over some line parallel to $x$ as $t$
runs over real numbers. This means that in the expansion of $I\circ\Phi(xt)$
all terms with nonzero powers of $t$ are proportional (parallel) to $x$.
We will write down some initial terms of this expansion dropping
the terms with zero power of $t$ and those obviously parallel to $x$:
$$I\circ\Phi(xt)=
\left(\frac{\Phi_3}{(x,x)}-\frac{2(\Phi_2,x)\Phi_2}{(x,x)^2}\right)t+$$
$$+\left(\frac{\Phi_4}{(x,x)}-\frac{2(\Phi_2,x)\Phi_3}{(x,x)^2}
-\frac{(\Phi_2,\Phi_2)\Phi_2+2(\Phi_3,x)\Phi_2}{(x,x)^2}+
\frac{4\Phi_2(x,\Phi_2)^2}{(x,x)^3}\right)t^2+\cdots$$

The terms with $t$ and $t^2$ must be proportional to $x$. The
proportionality conditions are polynomial relations in $x$. If they
hold for sufficiently many $x$'s in general position, then they hold
everywhere.

The coefficient with $t$ is equal to
$$\frac{\Phi_3}{(x,x)}-\frac{2(\Phi_2,x)\Phi_2}{(x,x)^2}.$$
Therefore, the map $\Phi_3(x,x)-2(\Phi_2,x)\Phi_2$ is everywhere
proportional to $x$. In particular, the inner product of this
map with $x$ is identically zero on the asymptotic cone $\{(x,x)=0\}$.
This implies that $(\Phi_2,x)=0$ on $C$. Hence $(\Phi_2,x)$ is divisible by
$(x,x)$, and so the map
$$\Phi_3-\frac{2(\Phi_2,x)\Phi_2}{(x,x)}$$
is a polynomial proportional to $x$. By Lemma \ref{pr} this polynomial is 
divisible by $x$ in the class of polynomials. Therefore, $\Phi_3$ is a linear 
combination with polynomial coefficients of $\Phi_2$ and $x$. So it always 
lies in the linear span of $\Phi_2$ and $x$.
In particular, $(\Phi_3,x)=0$ on $C$.

The term with $t^2$ is 
$$\frac{\Phi_4}{(x,x)}-\frac{2(\Phi_2,x)\Phi_3}{(x,x)^2}
-\frac{(\Phi_2,\Phi_2)\Phi_2+2(\Phi_3,x)\Phi_2}{(x,x)^2}+
\frac{4\Phi_2(x,\Phi_2)^2}{(x,x)^3}.$$
Multiply this expression by $(x,x)^2$ and restrict it to the asymptotic 
cone. We obtain that $\Phi_2(\Phi_2,\Phi_2)$ is parallel to $x$ on $C$
(note that all other terms are zero on the asymptotic cone).
This means that either $\Phi_2$ is parallel to $x$ on $C$ or the 
coefficient is zero. In both cases we have $(\Phi_2,\Phi_2)=0$ on $C$.
$\Box$

{\sc Example.} Let us construct an example of transformation that takes all 
lines to circles and has the identical differential at 0. Pick up a point
$a\in\R^n$ and consider the composition of the mirror reflection
$$x\mapsto x-2\frac{(a,x)a}{(a,a)}$$  
with respect to the orthogonal complement to $a$ and the inversion 
$$x\mapsto a+\frac{(a,a)(x-a)}{(x-a,x-a)}$$
with center $a$ and radius $|a|$ (so that 0 is fixed). Denote the resulting 
local diffeomorphism by $T^a$. We have 
$$T^a(x)=\frac{(a,a)x+(x,x)a}{(a,a)+2(a,x)+(x,x)}=x+\frac{(x,x)a-2(a,x)x}{(a,a)}+\cdots.$$
In particular, the quadratic term of $T^a$ has the form 
$$T^a_2(x)=\frac{(x,x)a-2(a,x)x}{(a,a)}$$
which is obviously parallel to $x$ on the asymptotic cone.

Now let us return to the general situation: we have a local diffeomorphism
$\Phi$ which rounds a sufficiently big and sufficiently general 
collection $\Lc$ of lines passing through 0.
Denote by $\Sc$ the corresponding set of circles.

\begin{prop}
\label{paral}
Suppose that $\Phi_2$ is parallel to $x$ on the asymptotic cone. Then 
all the circles from $\Sc$ pass through another common point different 
from the origin.
\end{prop}

To prove this, we need 2 very simple algebraic lemmas.

\begin{lemma}
\label{w}
Assume that a linear map $\Lambda:\R^n\to\Lambda^2\R^n$ satisfies
the condition $\Lambda(x)\wedge x=0$ everywhere. Then there is a
vector $b\in\R^n$ such that $\Lambda(x)=b\wedge x$. 
\end{lemma}

\proof Introduce a coordinate system $(x_0,\dots,x_{n-1})$ in $\R^n$.
Let $\Lambda_{ij}(x)$ be the coordinates of $\Lambda(x)$ in the
standard basis of $\Lambda^2\R^n$. These are linear functions in $x$.
The condition $\Lambda\wedge x=0$ can be written in coordinates
as follows:
$$\Lambda_{ij}x_k+\Lambda_{jk}x_i+\Lambda_{ki}x_j=0.\eqno{(*)}$$
The above formula implies that $\Lambda_{ij}$ vanishes on the
subspace $x_i=x_j=0$. Therefore, $\Lambda_{ij}=b_{ij}x_j-c_{ij}x_i$
where $b_{ij}$ and $c_{ij}$ are some numbers. Substitute this
expression to $(*)$:
$$(b_{ij}x_j-c_{ij}x_i)x_k+(b_{jk}x_k-c_{jk}x_j)x_i+
(b_{ki}x_i-c_{ki}x_k)x_j=0.$$
Equating the coefficient with $x_ix_j$ to zero we obtain
$b_{ki}=c_{jk}$. This implies that:
\begin{itemize}
\item the coefficient $b_{ki}$ is independent of $i$, denote it by $b_k$;
\item the coefficient $c_{jk}$ is independent of $j$, denote it by $c_k$;
\item $b_k=c_k$.
\end{itemize}
Now we have $\Lambda_{ij}=b_ix_j-b_jx_i$ which means that
$\Lambda(x)=b\wedge x$ where $b$ is the vector with coordinates
$(b_0,\dots,b_{n-1})$. $\Box$

Recall that a map $\Gamma:\C^n\to\C^n$ is {\em defined over reals}
if it takes $\R^n\subset\C^n$ to $\R^n$.

\begin{lemma}
\label{G}
Let $\Gamma:\C^n\to\C^n$ be a vector-valued quadratic form
(i.e., a homogeneous polynomial map of second degree) defined over reals
and such that $\Gamma(x)$ is everywhere parallel to $x$ on $C$.
Then $\Gamma$ has the form $\Gamma(x)=b(x,x)+\lambda(x)x$ where
$b\in\R^n$ and $\lambda$ is a linear functional. 
\end{lemma}

\proof 
Since $\Gamma$ is everywhere parallel to $x$ on the cone $C$, we have
$\Gamma(x)\wedge x=0$ there. Therefore, $\Gamma\wedge x$
is divisible by $(x,x)$. Denote the quotient by $\Lambda$.
It is a linear map from $\R^n$ to $\Lambda^2\R^n$.
Moreover, we have $\Lambda\wedge x=0$ because
$(\Gamma\wedge x)\wedge x=0$. By Lemma \ref{w} it follows that
$\Lambda=b\wedge x$ and hence $(\Gamma-b(x,x))\wedge x$ vanishes
everywhere. This means that the polynomial map $\Gamma-b(x,x)$
is proportional to $x$. By Lemma \ref{pr} we have
$\Gamma-b(x,x)=\lambda(x)x$ where $\lambda$ is some linear function.
$\Box$

{\sc Proof of Proposition \ref{paral}.}
By Lemma \ref{G} the second-order part $\Phi_2$ of a rectifying
diffeomorphism $\Phi$ has the form $\Phi_2(x)=b(x,x)+\lambda(x)x$
where $b$ is some vector from $\R^n$ and $\lambda$ is a linear
functional.

Consider a circle from $\Sc$ with the tangent vector $x$ at 0.
The acceleration with respect to the natural parameter is
$$2\frac{\Phi_2-\frac{(\Phi_2,x)x}{(x,x)}}{(x,x)}=
2\frac{\Phi_2-\lambda(x)x-(b,x)x}{(x,x)}=2\left(b-\frac{(b,x)x}{(x,x)}\right)$$
that is the same as for the circle passing through $b/(b,b)$.
But the circle is determined by its velocity $x/|x|$ and acceleration 
(both with respect to the natural parameter). It follows that all the circles
from $\Sc$ pass through $b/(b,b)$. $\Box$

Now we can give a simple proof of Theorem \ref{rect23}.

{\sc Proof of Theorem \ref{rect23}.} In dimensions 2 and 3 the asymptotic 
cone does not contain any plane. Therefore, $\Phi_2$ must be parallel to $x$
everywhere on the cone. Now Proposition \ref{paral} is applicable. $\Box$

{\sc Example.} In dimension 4, the statement of Theorem \ref{rect23} does not hold.
To construct a counterexample, introduce a complex structure on $\R^4$
and identify $\R^4$ with $\C^2$ by means of this complex structure.
Consider any complex projective transformation $\Phi$ preserving the origin.
It takes complex lines to complex lines, and on each line it induces
a projective transformation. On the other hand, a complex projective
transformation of a complex line takes real lines to circles.
Therefore $\Phi$ takes real lines to circles (note that each real
line belongs to exactly one complex line). Thus we get a rectifiable
family of circles (through 0). But these circles do not pass through a common
point different from the origin since different complex lines meet
only at the origin.

Theorem \ref{rect23} fails in dimension 4 by the following simple reason. 
The asymptotic cone now contains many planes, so there is no reason
anymore for $\Phi_2(x)$ to be everywhere parallel to $x$ on $C$. 
 
\section{Algebraic criteria for rectification}

We are going to prove now that the conditions on $\Phi_2$ stated in
Proposition \ref{cond} are not only necessary but also sufficient in a sense.

\begin{prop}
\label{suff}
If a vector-valued quadratic form $\Gamma:\C^n\to\C^n$ defined over reals
satisfies the conditions $(x,\Gamma(x))=(\Gamma(x),\Gamma(x))=0$ on the
asymptotic cone, then there exists a germ of diffeomorphism
$\Phi:(\R^n,0)\to(\R^n,0)$
that rounds lines passing through the origin and such that
$d_0\Phi=id$, $\Phi_2=\Gamma$, i.e., $\Phi=x+\Gamma$ up to third-order
terms.
\end{prop}

\proof Let us introduce the following notation:
$$\lambda=\frac{(\Gamma,x)}{(x,x)},\quad\mu=\frac{(\Gamma,\Gamma)}{(x,x)}.$$
We know that $\lambda$ and $\mu$ are polynomials in $x$ ($\lambda$ is a linear functional
and $\mu$ is a quadratic form). 

First assume that $\lambda=0$ (i.e., $(\Gamma,x)=0$ everywhere). Let us look
for a diffeomorphism $\Phi$ of 
the form $\Phi(x)=x+\Gamma(x)f(x)$ where $f$ is some smooth function that 
is equal to 1 at 0. We want $\Phi$ to take all lines (passing through 0) to 
circles. Denote by $I$ the inversion with center at 0 and 
radius 1. Then the germ of diffeomorphism 
$$I\circ\Phi=\frac{x+\Gamma f}{(x+\Gamma f,x+\Gamma f)}=
\frac 1{(x,x)}\frac{x+\Gamma f}{1+\mu f^2}$$ 
sends a neighborhood of 0 to a neighborhood of $\infty$ and is supposed to
take each line (passing through 0) to a parallel line. For that it suffices
to require that $f/(1+\mu f^2)=1$. Indeed, under the latter requirement
we have
$$I\circ\Phi(xt)=t^{-1}\frac x{(x,x)(1+\mu f(xt)^2)}+\frac{\Gamma}{(x,x)},$$
and the right-hand side has the form ``something parallel to $x$ plus
a term independent of $t$'' which means that $I\circ\Phi(xt)$ runs
over a line parallel to $x$ as $t$ runs over reals.
Solving the corresponding quadratic equation on $f$, we obtain 
$$f=\frac{1-\sqrt{1-4\mu}}{2\mu}.$$
We see that $f$ is a smooth analytic function near 0 such that $f(0)=1$
as we wanted.

Now suppose that $\lambda\ne 0$. Let us look for a diffeomorphism $\Phi$
of the form $\Phi=T^a\circ\Psi$ where $\Psi$ is
some other local diffeomorphism at 0. If $\Psi$ takes all lines passing
through 0 to circles, then the same is true for $\Phi$. We will try to kill
$\lambda$ by choosing an appropriate center $a$. For the second-order terms
we have $\Phi_2=\Psi_2+T^a_2$. So it suffices to take $a$ such that
$\lambda(x)=-(a,x)/(a,a)$.
Now $(\Psi_2,x)=0$ everywhere, so we reduced our problem to the previous case 
($\lambda=0$) which is done. $\Box$

Consider a simple bundle $\Sc$ of circles passing through 0 such that in each 
direction there goes a unique circle from $\Sc$. Such bundle 
is called {\em complete}. Now we can give a description of 
complete rectifiable bundles of circles in pure algebraic terms.

\begin{thm}
\label{class}
Complete rectifiable bundles of circles in $\R^n$ are in one-to-one
correspondence with quadratic homogeneous maps $\Gamma:\C^n\to\C^n$
defined over reals and satisfying the conditions
$(x,\Gamma(x))=(\Gamma(x),\Gamma(x))=0$ on the asymptotic cone, modulo
maps of the form $x\mapsto \lambda(x)x$ where $\lambda$ are linear
functionals.
\end{thm}

\proof To each complete rectifiable bundle $\Sc$ of circles assign the
quadratic part $\Phi_2$ of any rectifying diffeomorphism $\Phi$.
We know that any quadratic homogeneous map $\Phi_2$ defined over reals
and satisfying Proposition \ref{cond} can be obtained in this way.
Let us see to what extend the quadratic map $\Phi_2$ is unique.
We saw already that for each circle from $\Sc$ it is enough to know the 
acceleration at 0 with respect to the natural parameter. The acceleration
of the circle with the tangent vector $x$ is equal to
$$w(\Phi_2)=2\frac{\Phi_2-\frac{(\Phi_2,x)x}{(x,x)}}{(x,x)}.$$ 
But the above expression does not determine $\Phi_2$. It is easy to see that
if $\Phi_2$ and $\Phi'_2$ differ by $\lambda(x)x$ where $\lambda$ is a linear 
functional, then $w(\Phi_2)=w(\Phi'_2)$ so the corresponding families are
the same.
Indeed, it follows from the observation that $\Phi_2-(\Phi_2,x)x/(x,x)$ is
just the projection of $\Phi_2$ to the orthogonal complement of $x$. 
Vice versa, if $w(\Phi_2)=w(\Phi'_2)$, then $\Phi_2-\Phi'_2$ is everywhere
parallel to $x$ (since the projections to the orthogonal complement coincide). Hence
$\Phi_2-\Phi'_2=\lambda(x)x$ where $\lambda$ is a linear functional. 
$\Box$

{\sc Example.}  In dimension 4, the condition $(x,\Gamma(x))=
(\Gamma(x),\Gamma(x))=0$ on $C$ can be interpreted in terms
of algebraic geometry as follows. Denote
by $Q$ the projectivization of the asymptotic cone. This is a
nondegenerate quadratic surface in $\CP^3$. Each point of $Q$
belongs to 2 straight lines lying entirely in $Q$.

To describe all lines in $Q$ it is convenient to identify $Q$
with the image of the Segre embedding
$$\CP^1\times\CP^1\to\CP^3,\quad  ([u_0:u_1],[v_0:v_1])\mapsto
[u_0v_0:u_0v_1:v_0u_1:u_1v_1]$$
(recall that any nondegenerate quadratic surface in $\CP^3$ can
be mapped to any other by a complex projective transformation).
Under this embedding, each horizontal line $\CP^1\times\{p\}$
and each vertical line $\{p\}\times\CP^1$ get mapped to straight lines.
Hence we have 2 families of lines in $Q$ such that every point
of $Q$ belongs to a unique line from each family. These families
of lines are called {\em generating families of lines}.
For each generating family of lines in $Q$ there is the
corresponding {\em generating family of planes} in $C$.
So the cone $C$ is covered by 2 generating families of planes,
and every line in $C$ belongs to exactly one plane from each
generating family. 

The conditions $(x,\Gamma(x))=(\Gamma(x),\Gamma(x))=0$
on the asymptotic cone are equivalent to the following statement: the
subspace spanned by $x$ and $\Gamma(x)$ lies entirely in $C$. This means
that $\Gamma$ takes $x$ to another point of some line or plane
lying entirely in $C$. The map $\Gamma$ is homogeneous. Therefore, it
gives rise to a map from $\gamma:\CP^3\to\CP^3$ preserving the
projectivization $Q$ of the asymptotic cone $C$. We know that
for each point $q\in Q$ there is a line lying entirely in $Q$
and containing both $q$ and its image $\gamma(q)$. We will deduce from this
that $\Gamma$ preserves at least one of the generating families of lines in
$Q$ (maybe both), i.e., takes each line from some generating family to
itself. Indeed, being an algebraic map, $\gamma$ cannot ``switch'' from one
generating family to the other. Below is a formal proof of this statement.

\begin{lemma}
\label{gen}
The map $\gamma$ preserves at least one generating family of lines in $Q$.
\end{lemma}

\proof The surface $Q$ is isomorphic to $\CP^1\times\CP^1$ via the
Segre map. Hence $\gamma$ can be given by 2 algebraic maps
$$X:(x,y)\in\CP^1\times\CP^1\mapsto X(x,y)\in\CP^1,$$
$$Y:(x,y)\in\CP^1\times\CP^1\mapsto Y(x,y)\in\CP^1.$$
We know that for each point $(x,y)\in\CP^1\times\CP^1$ we have
$X(x,y)=x$ or $Y(x,y)=y$. Therefore, $Q$ is the union of 2 algebraic
subsets defined by the equations $X(x,y)=x$ and $Y(x,y)=y$.
Since $Q$ is irreducible, at least one of our equations is
satisfied identically, which means that $\gamma$ preserves at least
one of the generating families of lines in $Q$. $\Box$

Now we can deduce the following:

\begin{prop}
Polynomial homogeneous maps $\Gamma:\C^4\to\C^4$ satisfying the
conditions $(x,\Gamma(x))=(\Gamma(x),\Gamma(x))=0$ on the asymptotic cone
preserve some generating family of planes in $C$. 
\end{prop}

\section{Complex and quaternionic structures}

From now on we will work in 4-dimensional space $\R^4$. 
This section reviews not only well-known classical facts about complex
and quaternionic structures, but also their relation to the geometry of 
the asymptotic cone $C$. 

Recall that a {\em complex structure} in $\R^4$ is a linear operator
$I:\R^4\to\R^4$
such that $I^2=-1$. We will always assume that the complex structure $I$ is 
compatible with the Euclidean structure, i.e., $I$ preserves the
inner product.
A complex structure clearly defines an action of $\C$ on $\R^4$ via 
linear conformal maps. 
From the definition it follows immediately that $I$ must be skew-symmetric,
i.e., $(x,Iy)=-(Ix,y)$ for all $x,y\in\R^4$. In particular, $(Ix,x)=0$. 
Since the operator $I$ is defined over reals and $I^2=-1$, it should have
eigenvalues $i$ and $-i$, both with multiplicity 2.

Note that $I$ preserves the asymptotic cone $C$ (being an orthogonal
operator). In particular, all eigenvectors of $I$ belong to $C$. 
We know that $(Ix,x)=0$ everywhere and in particular on $C$.
From the conditions $(x,x)=(Ix,Ix)=(Ix,x)=0$ on $C$ it follows that
the subspace spanned by $x$ and $Ix$ lies entirely in $C$. 
Hence $I$ preserves one of the generating families of planes
in $C$.

On the other hand, the complex structure $I$ defines a canonical
orientation on $\R^4$. Let us recall the definition. Take 2 vectors
$x,y\in\R^4$ in general position. By definition, the canonical orientation
is the orientation of the basis $x,y,Ix,Iy$. This orientation is
well-defined (i.e., independent of the choice of $x$ and $y$) because the
set of degenerate pairs $(x,y)$ (such that $x,y,Ix,Iy$ are linearly
dependent) has real codimension 2 in the space $\R^8$ of all pairs. So we
can always avoid this set going from any nondegenerate pair to any other.
In fact, the degeneracy locus consists of all pairs $x,y$ that are linearly
dependent over $\C$, so it is a complex hypersurface. 

\begin{prop}
The space of all complex structures on $\R^4$ has 2 connected components. 
Complex structures from the same component preserve the same generating
family of planes in $C$ and provide the same canonical orientation. 
\end{prop}

A connected component to which a complex structure $I$ belongs will be called 
the {\em orientation} of $I$. Note that the orientation of $I$ has nothing to do with
$\det(I)$ which is always equal to 1 --- any complex structure preserves orientation
of the ambient space. 

Now let us pass to quaternionic structures. A {\em quaternionic structure} on
$\R^4$ is a choice of 3 linear operators $I,J,K:\R^4\to\R^4$ such that 
$$I^2=J^2=K^2=-1,$$
$$IJ=-JI=K,\quad JK=-KJ=I,\quad KI=-IK=J.$$
In particular, the operators $I,J,K$ are complex structures. We will assume
that they are compatible with the inner product. 
A quaternionic structure gives rise to an action of the skew-field 
$\H$ of quaternions on $\R^4$ via linear conformal maps. This action is
called the {\em quaternionic multiplication}. 

\begin{lemma}
Let $(I,J,K)$ be any quaternionic structure on $\R^4$. Then all 3 complex
structures $I,J,K$ have the same orientation. Therefore, quaternionic
multiplication preserves one of the generating families of planes in the
asymptotic cone.
\end{lemma}

\proof Let us prove, for example, that $I$ and $J$ provide the same canonical 
orientation. Take any vector $e\in\R^4$. It is enough to show that the 
bases $(e,Ke,Ie,IKe)$ and $(e,Ke,Je,JKe)$ have the same orientation. 
But $IKe=-Je$ and $JKe=Ie$, so the statement becomes obvious. $\Box$

Let $a\in\H$ be a quaternion. It gives rise to the operator of
multiplication $A:x\mapsto ax$. If $a=a_0+a_1i+a_2j+a_3k$, then the
corresponding operator is $A=a_0+a_1I+a_2J+a_3K$. We know that the operator
$A$ satisfies the conditions $(x,Ax)=(Ax,Ax)$ on $C$. In particular, both
forms $(Ax,Ax)$ and $(Ax,x)$ are divisible by $(x,x)$. We can write down
the quotients explicitly.

\begin{lemma}
\label{quot}
If $A$ is the operator of multiplication by a quaternion $a\in\H$ (with 
respect to some quaternionic structure on $\R^4$), then 
$$(Ax,Ax)=(a,a)(x,x),\quad (Ax,x)=\Re(a)(x,x).$$
In particular, these forms are independent of the choice of a quaternionic 
structure.  
\end{lemma}

\proof This is a very simple computation based on the fact that 
$(x,Ix)=(x,Jx)=(x,Kx)=0$ for all $x\in\R^4$. $\Box$

Let us summarize some properties of quaternionic structures that are of
particular importance for us. These properties follow directly from what
we saw already. 

\begin{prop}
\label{quatern}
The set of all quaternionic structures in $\R^4$ has 2 connected components.
Each component corresponds to a certain orientation of 3 complex structures
involved. Quaternionic multiplications with respect to quaternionic
structures from the same component preserve the same generating family of
planes in $C$.
Different components correspond to different families of generating planes.   
\end{prop}

We will say that quaternionic structures from the same connected component 
have the same {\em orientation}. Note that the orientation of a quaternionic 
structure has nothing to do with determinants of quaternionic multiplications.
Quaternionic multiplications (with respect to any quaternionic structure)
always preserve the orientation of the ambient space.

{\sc Example.} Identify $\R^4$ with $\H$. Denote by $I,J,K$ the operators of
left multiplication by $i,j,k$ respectively. The structure $(I,J,K)$ is 
called the {\em left quaternionic structure} on $\H$. If we take right
multiplication instead of left multiplication, then we get the {\em right
quaternionic structure}. Left and right quaternionic structures on $\H$
have different orientations.

Let us introduce some notions. We say that a linear operator is {\em almost 
orthogonal} if it has the form $const\cdot A$ where $A$ is an orthogonal operator.
Analogously, an operator is {\em almost skew-symmetric} if it has the form 
$const+A$ where $A$ is skew-symmetric. 

\begin{prop}
A linear operator $A:\R^4\to\R^4$ is the multiplication by a quaternion
(with respect to some quaternionic structure on $\R^4$) if and only if it
is almost orthogonal and almost skew-symmetric. 
The property of being a quaternionic multiplication depends only on the
orientation of a quaternionic structure, not on a structure itself.
\end{prop}

\proof A quaternionic multiplication is clearly almost orthogonal and almost 
skew-symmetric. If follows from Lemma \ref{quot}. Now consider an almost
orthogonal and almost skew-symmetric operator $A$ and present it by a
matrix in some orthonormal basis. Denote by $a_0,a_1,a_2,a_3$ the entries
of the first column of $A$. Since $A$ is almost skew-symmetric, it has the
form
$$\left(\begin{array}{cccc}
  a_0 & -a_1 & -a_2 & -a_3 \\
  a_1 & a_0  & \alpha & \beta \\
  a_2 & -\alpha & a_0 & \gamma \\
  a_3 & -\beta & -\gamma & a_0 
\end{array}\right).$$ 
The columns must be orthogonal and have the same length. From the
corresponding equations we obtain that either
$\alpha=a_3,\beta=-a_2,\gamma=a_1$ or $\alpha=-a_3,\beta=a_2,\gamma=-a_1$.
The first case corresponds to the left multiplication by
$a=a_0+a_1i+a_2j+a_3k$ with respect to the standard quaternionic
structure (assigned to the given basis). The second case corresponds to
the right multiplication by $a$. 
No matter what orthonormal basis we chose. Thus the second statement
follows. $\Box$

\section{Quaternionic fractional transformations}

Let us identify $\R^4$ with the skew-field $\H$ of quaternions. Consider 2
affine maps $A,B:\R^4\to\R^4$. The map $B^{-1}A$ (the multiplication and
the inverse are in the sense of quaternions) is called a (left)
{\em fractional quaternionic transformation} provided that it is one-to-one
at least in some open subset of $\R^4$. A {\em right quaternionic fractional
transformation} is a local transformation of the form $AB^{-1}$ where $A$
and $B$ are some affine maps.

{\sc Example 1.} Any real projective transformation is quaternionic
fractional. This corresponds to the case when $B$ takes real values only.

{\sc Example 2.} Any complex projective transformation is quaternionic fractional.
This happens if $B$ takes complex values only and both $A$ and $B$ are 
complex linear (i.e., commute with the multiplication by $i$). 

{\sc Example 3.} Consider a map of the form $x\mapsto (xa+b)^{-1}(xc+d)$
where $a,b,c,d$ are quaternions. We are assuming that the denominator 
is not proportional to the numerator (in particular, the denominator is not
identically zero).
This map is called a (left) {\em quaternionic projective transformation}.
Any quaternionic projective transformation is clearly quaternionic
fractional. Note that each quaternionic projective transformation takes all
lines to circles. Indeed, we have 
$$(xa+b)^{-1}(xc+d)=(xa+b)^{-1}((xa+b)\alpha+\beta)=\alpha+(xa+b)^{-1}\beta$$
where $\alpha=a^{-1}c$, $\beta=d-b\alpha$. Hence a quaternionic 
projective transformation is a composition of a dilatation, reflected
inversion and a translation. This composition obviously rounds lines. 
 
\begin{prop}
Any quaternionic fractional transformation rounds lines (to be more precise:
it takes germs of lines to germs of circles).  
\end{prop}

\proof Consider a line $L$ in $\R^4$. Let $t$ be a linear parameter on $L$.
If $A$ and $B$ are some affine maps, then their restrictions to $L$ 
are $at+b$ and $ct+d$ respectively. So on the line $L$ the transformation
$A^{-1}B$ coincides with the quaternionic projective transformation 
$x\mapsto (ax+b)^{-1}(cx+d)$. But the latter rounds lines. $\Box$  
 
\section{Rectification at a point}

In this section, we will prove the following theorem:

\begin{thm}
\label{rect4}
Consider a simple bundle of circles in $\R^4$ containing sufficiently many
circles in general position. If this bundle is
rectifiable, then there exists a left or right quaternionic fractional
transformation $T$ such that $T^{-1}$ sends all these circles to straight
lines.
\end{thm}

Denote the given set of circles by $\Sc$. Let $\Phi$ be
a local diffeomorphism such that $d_0\Phi=id$ and $\Phi^{-1}$ rectifies all 
circles from $\Sc$. Then by Proposition \ref{cond} the 
quadratic term $\Phi_2$ satisfies the relations
$(\Phi_2,x)=(\Phi_2,\Phi_2)=0$
on the asymptotic cone. This means that $\Phi_2$ preserves one of the 
generating families of planes in $C$.

\begin{lemma}
\label{Phi_2}
There exists a linear operator $A:\R^4\to\R^4$ such that 
$\Phi_2(x)=A(x)x$ or $\Phi_2(x)=xA(x)$ where the product is in the sense 
of quaternions.
\end{lemma}

\proof
Fix an identification $\R^4=\H$. Extend the operators $I$, $J$ and $K$ of
left multiplication by $i$, $j$ and $k$ respectively to $\C^4$ by complex
linearity. Note that the operator $I$ is quite different from the
multiplication by $\sqrt{-1}$ in $\C^4$. By Proposition \ref{quatern} the
left quaternionic multiplication preserves one of the generating families of
planes in $C$. Assume that $\Phi_2$ preserves the same family. Otherwise we
should consider the right multiplication instead of the left multiplication. 
 
Recall that the {\em quaternionic conjugation} is the map
$$x=x_0+x_1i+x_2j+x_3k\mapsto \bar x=x_0-x_1i-x_2j-x_3k.$$
We can extend this map to $\C^4$ by complex linearity. Note that 
$i$ is now a vector from $\R^4$, not a complex number.
Let us multiply $\Phi_2$ by $\bar x$ in the sense of quaternions.
Note that 
$$\Phi_2\bar x=(\Phi_2,x)+(\Phi_2,Ix)i+(\Phi_2,Jx)j+(\Phi_2,Kx)k.$$
But this expression is zero on the cone $C$ since $\Phi_2$, $x$, $Ix$,
$Jx$ and $Kx$ lie on the same plane belonging to $C$.
Therefore, $\Phi_2\bar x$ is divisible by $(x,x)$. The quotient is a linear
map $A$. Since $\bar x/(x,x)=x^{-1}$, we have $\Phi_2x^{-1}=A(x)$, i.e., 
$\Phi_2(x)=A(x)x$. $\Box$

Now we can prove Theorem \ref{rect4} and even more precise statement:

\begin{thm}
Under assumptions of Theorem \ref{rect4} the family of circles can be
obtained from the family of their tangent lines by one of the transformations  
$x\mapsto (1-A(x))^{-1}x$ or $x\mapsto x(1-A(x))^{-1}$ where $A$ is some 
linear operator. This answer does not depend on the choice of a quaternionic 
structure. 
\end{thm} 

\proof Note that both transformations have the identical derivative at 0
and their second-order terms are $A(x)x$ and $xA(x)$ respectively. 
These transformations are quaternionic fractional so they round lines. 
The corresponding families of circles passing through 0 are determined by 
the second-order terms. But by Lemma \ref{Phi_2} the quadratic maps $A(x)x$
and $xA(x)$ are the only possible second-order terms of transformations
that round lines. $\Box$
 
For a complete rectifiable bundle $\Sc$ of circles there is a transformation 
of the form  
$x\mapsto (1-A(x))^{-1}x$ or $x\mapsto x(1-A(x))^{-1}$ that takes
the family of all lines passing through 0 to $\Sc$. To fix the idea,
assume that this is the left transformation $\Phi:x\mapsto (1-A(x))^{-1}x$.

\begin{prop}
The center of the circle from $\Sc$ with the tangent vector $x$ at 0 is 
$-\frac 12(\Im A(x))^{-1}x$. This point can be infinite which means that the 
corresponding circle is a straight line. 
\end{prop} 

\proof We know that the acceleration with respect to the natural parameter is 
$$w(x)=2\frac{\Phi_2-\frac{(\Phi_2,x)x}{(x,x)}}{(x,x)}.$$
Therefore the center is located in the point
$$\frac w{(w,w)}=
\frac 12\frac{ \frac{\Phi_2}{(x,x)} - \frac{(\Phi_2,x)x}{(x,x)^2} }
{ \frac{(\Phi_2,\Phi_2)}{(x,x)^2} - \frac{(\Phi_2,x)^2}{(x,x)^3} }.$$
By Lemma \ref{quot} we have $(\Phi_2,\Phi_2)=(A,A)(x,x)$ and
$(\Phi_2,x)=(\Re A)(x,x)$.
Simplifying the above expression we get the following formula for the center:
$$\frac 12\left(\frac{A-\Re A}{(A,A)-(\Re A)^2}\right)x=
\frac 12\frac{\Im(A)}{(\Im A,\Im A)}x=-\frac 12(\Im A)^{-1}x.$$
$\Box$

The previous proposition has the following geometric corollary:

\begin{cor}
The family $\Sc$ contains at least one line. The union of all straight
lines from $\Sc$ is a vector subspace of $\R^4$.
\end{cor}

{\sc Remark.} We see that the set of all complete rectifiable families of
circles passing through 0 is naturally identified with the union of 2 affine 
spaces of dimension 12 (=dimension of all possible $\Im A(x)$). 
The intersection of these components has dimension 4 and consists of all
families rectifiable by means of inversions (i.e., of 
families whose circles meet at a point different from 0 --- this happens
if $\Im A$ is independent of $x$).
The two components can be distinguished by the ``orientation''.

We can describe an affine structure on each component in geometric terms.
Namely, take any 2 circles $S_1$ and $S_2$ tangent at 0. After an inversion, 
they become parallel lines. For two parallel lines $L_1$ and $L_2$ we can
take their barycentric combination
$$L=\lambda L_1+(1-\lambda)L_2=
\{\lambda x+(1-\lambda)y|\ x\in L_1,y\in L_2\},\quad \lambda\in\R.$$
Make the inversion again. The line $L$ goes to a circle $S$. Put by
definition $S=\lambda S_1+(1-\lambda)S_2$. Now we can take barycentric
combinations of complete bundles of circles. Namely, let the circle of the
new bundle passing through 0 in direction $x$ be
$S=\lambda S_1+(1-\lambda)S_2$ where $S_1$ and $S_2$ are
circles from the old bundles going from 0 in direction $x$.    
It turns out that if two rectifiable bundles have 
the same ``orientation'', then their barycentric combinations are also
rectifiable.
   
{\sc Open question.}
How many complete rectifiable simple bundles of circles are there? We saw
that in $\R^n$ the space of all complete rectifiable bundles of
circles passing through 0 is finite-dimensional. What is its dimension
(as a function of $n$)? Is there an explicit geometric description of such
bundles in dimensions $>4$?

\end{document}